\documentclass[11pt]{amsart}
\usepackage{enumitem, verbatim}
\usepackage{amsmath, amssymb, amsthm}

\usepackage[margin=2 cm,heightrounded=true,centering]{geometry}

\newtheorem{theorem}{Theorem}[section]
\newtheorem{proposition}[theorem]{Proposition}
\newtheorem{lemma}[theorem]{Lemma}
\newtheorem{remark}[theorem]{Remark}
\newtheorem{corollary}[theorem]{Corollary}
\newtheorem{definition}[theorem]{Definition}

\newcommand{\abs}[1]{\left|#1\right|}

\newcommand{\metr}[1]{\langle #1 \rangle}
\newcommand{\pdx}[2][]{\frac{\partial #1}{\partial #2}}
\newcommand{\ddt}[1][]{\frac{d #1}{dt}}
\newcommand{\brax}[1]{\left( #1 \right)}

\newcommand{\calL}{\mathcal{L}}

\newcommand{\at}{\Big|}

\renewcommand{\tilde}{\widetilde}

\DeclareMathOperator{\ric}{Ric}
\DeclareMathOperator{\riem}{Rm}
\DeclareMathOperator{\scc}{Sec}
\DeclareMathOperator{\hess}{Hess}

\DeclareMathOperator{\interior}{int}
\DeclareMathOperator{\tr}{tr}
\DeclareMathOperator{\tf}{tf}
\DeclareMathOperator{\dist}{dist}
\DeclareMathOperator{\divergence}{div}
\DeclareMathOperator{\diag}{diag}

\begin{document}

\title{Bakry-\'Emery Ricci Curvature Bounds on Manifolds with Boundary}
\author{Kenneth Moore}
\address{Department of Mathematical and Statistical Sciences,
University of Alberta, Edmonton, Alberta, Canada T6G 2G1}
\email{kjmoore1@ualberta.ca}
\author{Eric Woolgar}
\address{Department of Mathematical and Statistical Sciences and Theoretical Physics Institute,
University of Alberta, Edmonton, Alberta, Canada T6G 2G1}
\email{ewoolgar@ualberta.ca}

\begin{abstract}
We prove a Bakry-\'Emery generalization of a theorem of Petersen and Wilhelm, itself a generalization of a theorem of Frankel, that closed minimal hypersurfaces in a complete manifold with a suitable curvature bound must intersect. We then prove splitting theorems of Croke-Kleiner type for manifolds bounded by hypersurfaces obeying Bakry-\'Emery curvature bounds. Motivated in part by the near-horizon geometry programme of general relativity, we do not assume that the Bakry-\'Emery vector field is of gradient type.
\end{abstract}

\maketitle

\section{Introduction}
\setcounter{equation}{0}

\subsection{Preliminaries} Frankel's theorem \cite{Frankel1} implies that compact totally geodesic submanifolds in a manifold with positive sectional curvature must intersect if the sum of the dimensions of the submanifolds equals or exceeds that of the ambient manifold. It was subsequently shown \cite{Frankel2, PW} that any two minimal hypersurfaces in a manifold of positive Ricci curvature must always intersect. In the case that the Ricci curvature is only nonnegative, the minimal surfaces need not intersect, but this can occur only in special situations. Rigidity theorems for the ambient manifold in this circumstance were obtained in \cite{Ichida}, \cite{Kasue}, \cite{Galloway}, \cite{CK}, and \cite{PW}. These results place topological restrictions on manifolds with nonnegative Ricci curvature that admit a minimal hypersurface. A version for the case of gradient Ricci solitons can be found in \cite[Theorem 7.4]{WW}. In this paper we prove similar results under general Bakry-\'Emery Ricci curvature and mean curvature assumptions, without requiring the Bakry-\'Emery 1-form $X$ to be exact.

The study of Bakry-\'Emery Ricci curvature was presaged by the study of manifolds with a positive density function $e^{-f}$, which goes back at least to Lichn\'erowicz \cite{Lichnerowicz}. Interest has revived in recent years, due to connections to a number of timely topics including Ricci solitons, general relativity, and smooth metric measure spaces.
Many results that apply to metrics with Ricci curvature bounds have been generalized to metrics with Bakry-\'Emery Ricci bounds, such as Myers's theorem \cite{Qian} and the Cheeger-Gromoll splitting theorem. Versions of the splitting theorem have been proved in \cite{FLZ, Wylie, KWW, Sakurai1, Sakurai2, Sakurai3, Lim1}. The Bakry-\'Emery Ricci tensor is defined as follows.

\begin{definition}\label{definition1.1}
Let $X$ be a vector field on a Riemannian manifold $(M^n,g)$. To minimize the burden of notation, we also use $X$ to denote the metric-dual $1$-form. The \emph{$m$-Bakry-\'Emery Ricci curvature} is
\begin{equation}
\label{eq1.1}
\begin{split}
\ric_X^m:=&\, \ric + \frac{1}{2}\pounds_Xg -\frac{1}{m}X\otimes X\ ,\ m\neq 0\ ,\\
\ric^0:=&\, \ric\ ,\ m=0\text{ and } X\equiv 0\ ,\\
\ric_X^{\infty}:=&\, \ric + \frac{1}{2}\pounds_Xg
\end{split}
\end{equation}
where $\pounds_Xg$ is the Lie derivative of $g$ along $X$. If $X=\nabla f$ for $f:M\to {\mathbb R}$, it is conventional to write $\ric_f^m$ or $\ric_f^{\infty}$, respectively, and then we take $f$ to be constant in the $m=0$ case.
\end{definition}

There is also a notion of Bakry-\'Emery mean curvature for hypersurfaces. Let $\Sigma$ be a hypersurface with a well-defined unit normal field. Let $\nu$ be any smooth extension of that field to an open neighbourhood about $\Sigma$ (the following definition does not depend on the choice of smooth extension). We choose signs so that the mean curvature $H$ is the trace of the shape operator; i.e., the trace of the map $S$ that sends $U\in T_p \Sigma$, $p\in \Sigma$, to $\nabla_U \nu\in T_p\Sigma$.
\begin{definition}\label{definition1.2}
The $X$-mean curvature of the hypersurface $\Sigma$ with respect to the (extended) unit normal field $\nu$ is defined to be
\begin{equation}
\label{eq1.2}
H_X=H-g(X,\nu)\ .
\end{equation}
If $H_X=0$ pointwise, we call $\Sigma$ a \emph{Bakry-\'Emery $X$-minimal hypersurface}. If $X=\nabla f$, we write $H_f=H-g(\nabla f, \nu)$.  A hypersurface with $H_f=0$ is called \emph{$f$-minimal}.
\end{definition}
When $\Sigma$ is a component of the boundary $\partial M$, we will choose the unit normal field to point into $M$. Our conventions are such that if $M$ is ${\mathbb R}^n$ minus a unit ball, the mean curvature of the unit $n$-sphere boundary is $(n-1)$.

\subsection{Results I. A Frankel-type theorem}

\begin{theorem}
\label{theorem1.3}
Let $M$ be an $n$-dimensional Riemannian manifold that admits an $m\in (0,\infty]$ and a vector field $X$ such that the Bakry-\'Emery Ricci curvature obeys $\ric_X^m > 0$. Let $N_1$ and $N_2$ be closed $X$-minimal hypersurfaces in $M$. Then $N_1$ and $N_2$ intersect.
\end{theorem}

Wei and Wylie \cite[Theorem 7.4]{WW} proved this theorem for the case of $f$-minimal hypersurfaces $H_f=0$ in complete manifolds with $\ric_f^{\infty}>0$ and $f$ bounded above. An alternative proof for this case was given in \cite{LW}.

Theorems of this type have an obvious topological consequence, for say that the manifold of Theorem \ref{theorem1.3} is nontrivially covered by another manifold such that the covering map is a local isometry and the covering manifold contains disjoint copies of the $X$-minimal hypersurface. Such a covering space  would violate the theorem, so coverings of this sort cannot occur. But if the curvature condition were slightly relaxed, becoming $\ric_X^m \ge 0$, then this possibility can arise. However, there remain topological restrictions. For what herein would be the $m=0$ case, see \cite[Theorem 4]{PW}. We give the $m>0$ generalization of this result in Corollary \ref{corollary4.2}. Galloway \cite{Galloway} and, more recently, Choe and Fraser obtain topological results in the case of a single minimal surface \cite[Theorem 2.5]{CF} (with $m=0$), essentially by passing to a covering space where the analysis of \cite{PW} would apply. We generalize the statement of \cite[Theorem 2.5]{CF} to the Bakry-\'Emery $m>0$ setting in the following theorem:
\begin{theorem}
\label{theorem1.4}
Let $M$ be a compact manifold with a vector field $X$ and an $m\in (0,\infty]$ such that $\ric_X^m \geq 0$ pointwise on $M$, and such that $M$ has a closed embedded $2$-sided $X$-minimal hypersurface $N$.
\begin{itemize}
\item[a)]
If $N$ is nonseparating, then $M$ is isometric to a mapping torus
\begin{equation}
\label{eq1.3}
\frac{N\times [0,a]}{(x,0)\sim (\phi(x),a)}
\end{equation}
where $\phi:\Sigma\to \Sigma$ is an isometry, and if $m>0$ is finite then $X$ is tangent to $\Sigma$, which is then a minimal hypersurface in the usual sense.

\item[b)] If $\Sigma$ is separating, let $D_1$ and $D_2$ be the connected components of $M\backslash \Sigma$. Then for $j=1,2$, the maps
\begin{equation}
\label{eq1.4}
i_*:\pi_1(\Sigma)\to\pi_1(\overline{D_j})\text{, } I_*:\pi_1(D_j)\to \pi_1(M)\text{, and }{\mathcal I}_*:\pi_1(\Sigma)\to \pi_1(M)
\end{equation}
induced by inclusion are all surjective.
\end{itemize}
\end{theorem}

Choe and Fraser point out that a consequence of part (b) in the $3$-dimensional case is that $\Sigma$ is then a Heegaard surface, dividing $M$ into two handlebodies. Furthermore, they then show that if the ambient manifold is a topological $3$-sphere with $\ric\ge 0$ then $\Sigma$ is unknotted (i.e., deformable to a standard embedding) \cite[Corollary 2.6]{CF}. The same result follows from Theorem \ref{theorem1.4} if $\ric_X^m\ge 0$.

We will employ a direct method of proof of Theorem \ref{theorem1.4}. But in the case of part (a) the fundamental group splits as a semi-direct product $\pi_1(M)\simeq \pi_1(N)\rtimes {\mathbb Z}$. By the Bakry-\'Emery version of Myers's theorem, this cannot occur when $\ric_X^m\ge \lambda g$ for some $\lambda>0$, for then the fundamental group must be finite. Hence, under the assumptions of the theorem, part (a) can arise only if $\ric_X^m(Y,Y)=0$ for some $Y\in T_p M\backslash \{ 0 \}$, $p\in M$. This observation points toward a different method of proof from the one we will employ. If the manifold obeys $\ric\ge 0$ but does not have finite fundamental group, then often there is a rigidity which forces the manifold to split. If the manifold admits a line (an inextendible minimizing geodesic), this is the well-known Cheeger-Gromoll theorem \cite{CG}. In the case of manifolds admitting a compact hypersurface with zero mean curvature, a splitting theorem applies (see, e.g., \cite{CK}). We now describe versions of the splitting theorem in our setting.

\subsection{Results II. Splitting theorems}
Sakurai has given versions of the splitting theorem that apply given Bakry-\'Emery Ricci and mean curvature bounds in place of the usual Ricci and mean curvature bounds \cite{Sakurai1, Sakurai2, Sakurai3} when the Bakry-\'Emery vector field is gradient, i.e., $X=\nabla f$. A simple modification allows us to obtain splitting theorems for nongradient $X$. We obtain versions of these theorems for all $m\in (0,\infty]\cup (-\infty,1-n]$.


We first generalize the classic warped product splitting \cite[Theorem 1]{CK} (see \cite[Theorem, p 169]{Ichida} and \cite[Theorem B.(1)]{Kasue} as well). For gradient $X$ with zero Bakry-\'Emery curvature bounds and $m\in (0,\infty]$, see \cite[Theorem 6.13]{Sakurai2}. For $m\in (-\infty,1-n]$ (where $n=\dim M$), we treat only the case of zero curvature bounds, but for the gradient case including nonzero bounds see \cite[Theorem 5.10]{Sakurai3}.

\begin{theorem}\label{theorem1.5}
Let $M$ be a complete manifold-with-boundary, with boundary components $N_1$ and $N_2$, at least one of which is compact. Suppose that there is an $m\in (0,\infty)$ and an $X$ such that $\ric_X^m(M)\geq -(n-1)\delta$ for $\delta\in \{ 0,1\}$. Suppose that the Bakry-\'Emery mean curvature of $N_1$ is $\leq - \sqrt{(n-1)(n+m-1)}\delta$ and of $N_2$ is $\leq \sqrt{(n-1)(n+m-1)}\delta$. Then $M$ is isometric to $N_1\times [0,\ell]$ with the metric $ds^2=dt^2+e^{2c\delta t}g_1$ where $g_1$ is the metric on $N_1$, $c:=\sqrt{\frac{n-1}{n+m-1}}$, and $\ric_{X^{\sharp}}^m(g_1)\ge -\frac{(n-1)^2}{(n+m-1)}\delta$ where $X^{\sharp}$ denotes the restriction of $X$ to $TN_1$. For $\delta=0$, the splitting is a Riemannian product and the projection of $X$ along $\partial_t$ vanishes. For $\delta=1$, the splitting is that of a warped product and the projection of $X$ along $\partial_t$ is constant, namely $g(X,\partial_t)=m\sqrt{\frac{n-1}{n+m-1}}$.
\end{theorem}

\begin{theorem}\label{theorem1.6}
Let $M$ be a complete manifold-with-boundary, with boundary components $N_1$ and $N_2$, at least one of which is compact. Suppose that $\ric_X^{m}(M)\geq 0$ where $m\in\{ \infty \}\cup (-\infty, 1-n]$. Further suppose that the Bakry-\'Emery mean curvatures of $N_1$ and of $N_2$ are each $\leq 0$. Then $M$ is isometric to $N_1\times [0,\ell]$ with the metric $ds^2=dt^2+g_1$ where $g_1$ is the metric on $N_1$, and $\ric_{X^{\sharp}}^m(g_1)\ge 0$. If $m\neq 1-n$ the projection of $X$ along $\frac{\partial}{\partial t}$ vanishes.
\end{theorem}

Next we generalize the half-space rigidity theorem (see \cite[Theorem 2]{CK}). Other authors have also proved Bakry-\'Emery versions \cite{Sakurai1, Sakurai2, Sakurai3, Yang}, under various different assumptions.

\begin{theorem}\label{theorem1.7}
Let $M$ be a complete manifold with non-empty compact connected boundary $N=\partial M$ and an asymptotic end. Suppose that there is an $m\in (0,\infty)$ and an $X$ such that $\ric_X^m(M)\geq -(n-1)\delta$ for $\delta\in \{ 0,1\}$. Suppose that the Bakry-\'Emery mean curvature of $N$ is $H_X \leq - \sqrt{(n-1)(n+m-1)}\delta$. Then $M$ is isometric to $[0,\infty)\times N$ with the metric $ds^2=dt^2+e^{2c\delta t}g_N$ where $g_N$ is the metric on $N$ and $c:=\sqrt{\frac{n-1}{n+m-1}}$. For $\delta=0$, the splitting is a Riemannian product and the projection of $X$ along $\partial_t$ vanishes. For $\delta=1$, the splitting is that of a warped product and the projection of $X$ along $\partial_t$ is constant, namely $g(X,\partial_t)=m\sqrt{\frac{n-1}{n+m-1}}$.
\end{theorem}

Again, for $m\in (-\infty,1-n]$ (where $n=\dim M$), we treat only the case of zero curvature bounds.

\begin{theorem}\label{theorem1.8}
Let $M$ be a complete manifold with non-empty compact connected boundary $N=\partial M$ and an asymptotic end. Suppose that there is an $m\in (-\infty,1-n]\cup \{ \infty\}$ and an $X$ such that $\ric_X^m(M)\geq 0$. Suppose that the Bakry-\'Emery mean curvature of $N$ is $H_X \leq 0$. Further suppose that $\left \vert \int\gamma X\cdot ds\right \vert\le K$ for some constant $K>0$ and every geodesic $\gamma$.
\begin{itemize}
\item [a)] If $m\neq 1-n$ then $M$ is isometric to $[0,\infty)\times N$ with the metric $ds^2=dt^2+g_N$ where $g_N$ is the metric on the leaves diffeomorphic to $N$ and is independent of $t$, and $\ric_{X^{\sharp}}^m(g_N)\ge 0$ where $X^{\sharp}$ denotes the restriction of $X$ to $TN$. The projection of $X$ along $\partial_t$ vanishes.
\item [b)] If $m=1-n$ then $M$ is isometric to $[0,\infty)\times N$ with the twisted product metric \begin{equation}
\label{eq1.5}
ds^2=dt^2+e^{\frac{2}{(n-1)}\int_0^t g(\gamma',X)ds}g_N.
\end{equation}
If $X$, considered as a 1-form, is closed, then $g(\gamma',X)$ is constant on the leaves $t=const$ and \eqref{eq1.5} is a warped product.
\end{itemize}
\end{theorem}

\subsection{Motivation} Consider the \emph{vacuum near horizon geometry (NHG) problem}. The problem is to classify solutions of $\ric_X^m=\lambda g$ for $\lambda\in{\mathbb R}$ and $m=2$ on closed manifolds $M$. If $X$, regarded as a 1-form, is closed, this is the \emph{static} vacuum NHG problem. If $X$ is closed and $\lambda\ge 0$, then it is known that $X$ vanishes and the problem reduces to solving the Einstein equation on $M$. Solutions of the vacuum NHG problem are candidate Killing horizons for static (if $X$ is closed) and stationary (if $X$ is not closed) degenerate (also called extreme) black holes in general relativity.

For the $\lambda=0$ static vacuum NHG problem, we can restrict attention to compact Einstein $n$-manifolds \cite{CRT}, so the results of \cite{CF} apply. Take $n<8$ and assume that the homology group $H_{n-1}(M)$ is non-trivial. Then each non-trivial class is represented by a smooth minimal hypersurface $\Sigma$. If this surface is nonseparating, then by \cite[Theorem 2.5]{CF}, $M$ is isometric to a mapping torus. The fundamental group is then $\pi_1(M)\cong \pi_1(\Sigma)\rtimes{\mathbb Z}$. If instead $\Sigma$ is separating, then $\pi_1(\Sigma)\to\pi_1(M)$ is surjective. Compare to \cite[Theorem 1]{KWW} (which does not have a homology assumption) when $\pi_1(\Sigma)$ is finite.

With this in mind, we ask whether these methods can shed light on the NHG problem when $X$ is not exact. This arises both in the non-static case and in the static $\lambda<0$ case. Unfortunately, the authors are not aware of general arguments that imply the existence of $X$-minimal surfaces except when $X$, considered as a 1-form, is exact: $X=df$. Thus one may view Theorem \ref{theorem1.4} as giving criteria for the non-existence of $X$-minimal surfaces. If on the other hand the existence of such surfaces could separately be established under reasonable conditions, then a rich variety of results would likely follow, mimicking results that follow from the interplay between existence results for ordinary minimal surfaces and ordinary Ricci curvature bounds. It is hoped that if a general existence theory is not available, a limited theory may arise within the context of the most relevant applications. Similar considerations apply to splitting theorems, where the Bakry-\'Emery Ricci bound may be nonzero and the relevant hypersurfaces are then of constant nonzero $X$-mean curvature.

A motivating example is the metric $g=g_{\Sigma}+dz^2$ on $\Sigma\times S^1$, with $\Sigma$ a compact hyperbolic $2$-manifold with metric $g_{\Sigma}$ of constant sectional curvature $-1$. Here $z$ is a coordinate on (a patch for) $S^1$. This metric obviously splits as a product, but would seem not to arise as a rigidity case in the Croke-Kleiner theorem \cite{CK} because the Ricci tensor is $\ric=\diag(-1,-1,0)\ge -1$ while the surfaces $\Sigma$ have mean curvature $0$, so the mean curvature and Ricci bounds do not match. To see that this yields a solution of the NHG problem \cite{Lim2}, choose the $1$-form $X=\sqrt{m}dz$. Then $X$ extends to be globally defined, though the coordinate $z$ is not, and so $X$ is closed but not exact. The Bakry-\'Emery Ricci tensor is $\ric_X^m=\diag(-1,-1,-1)=-g$, so $\ric_X^m$ has equal eigenvalues. Setting $m=2$, we have solved the static vacuum NHG problem with $\lambda=-1$ although $g$ is not Einstein. The level sets $z=const$ are not $X$-minimal, but are constant $X$-mean curvature surfaces with $H_X=\pm 1$ depending whether the normal direction is chosen parallel or anti-parallel to $X$. This manifold obeys all the conditions of a splitting theorem that we will establish herein.

In our example, we can now deduce that there are no pairs $(g',X')$ that obey $\ric_{X'}^m\ge -g'$ and that differ from $(g,X)$ only on a region of diameter less than the period of the $S^1$ factor (so that at least one $X$-minimal surface remains intact). If analytical techniques were available that could ensure that a constant $X$-mean curvature surface were always present after a general variation with no diameter restriction or other restrictive assumption, we would obtain a stronger result. Herein we ask only what geometrical properties follow when such surfaces are present.

There are claims in the literature that the static vacuum NHG problem has no non-Einstein solutions when $\ric_X^{m=2}=\lambda g$, $dX=0$, $\lambda\le 0$ (see \cite[Theorem 4.1]{KL}). As the above example shows \cite{Lim2}, such claims are not correct. The error can be traced to \cite{CRT} which establishes the $\lambda=0$ case and states that the reasoning also holds when $\lambda<0$. But in fact the reasoning fails when $\lambda<0$ if $\divergence X -|X|^2-2\lambda=0$ because in this case the condition $dX=0$ is not sufficient to ensure that $X$ is exact. Thus, the static vacuum near horizon geometry problem with $\lambda<0$ remains open, as does the general NHG problem for non-closed $X$. The NHG problem is itself an attempt to understand a larger open problem, that of black hole uniqueness and non-uniqueness in spacetime dimensions greater than $4$.

\subsection{Organization of the paper}
Section 2.1 contains the proof of Theorem \ref{theorem1.3}. Section 2.2 contains 
a topological result
and a generalization of Theorem \ref{theorem1.3} when the curvature bounds are not zero. Section 3.1 describes a Bochner identity on functions, while Section 3.2 describes the Riccati (Raychaudhuri) equation with Bakry-\'Emery curvature bounds. Section 3.3 uses the Riccati equation for $X$-mean curvature to derive comparison-type lemmata which are used in Sections 4 and 5. All of these results are established without requiring the Bakry-\'Emery vector field $X$ to be gradient. This is a feature of the Riccati equation approach---see especially Lemma \ref{lemma3.3}---which uses only line integrals defined along a chosen geodesic and so does not need $X$ to be globally integrable. Section 4 uses the analysis of Section 3 to give the proof of Theorem \ref{theorem1.4} and a generalization of \cite[Theorem 4]{PW}. When $X$ is gradient, these results can separately be obtained from splitting theorems. We prove the splitting theorems without the gradient assumption in Section 5.

\subsection{Conventions} Since this paper generalizes known results in the $m=0$ case which would often require special handling, we omit that case throughout. We define the Laplacian $\Delta$ to be the trace of the Hessian. We use $\tr T$ and $\tf T$ to denote the trace and the tracefree part, respectively, of a $(1,1)$-tensor $T$ (if $T$ is, say, a $(0,2)$-tensor, then it is implicit that we raise an index with $g$ before tracing). We define mean curvature $H$ of a given hypersurface with respect to a unit normal field $\nu$ to be $H=\divergence \nu = \tr (\nabla \nu)$, so that the mean curvature of a round sphere in ${\mathbb R}^3$ with respect to the normal field pointing to infinity is positive. In particular, if $d$ is a distance function with a smooth level set $d=c$, its mean curvature will be $H=\divergence \nabla d=\Delta d$ evaluated at $d=c$.

\subsection{Acknowledgements} This research was funded by NSERC Discovery Grant RGPIN–2017–04896 to EW. We are grateful to GJ Galloway for bringing references \cite{Ichida}, \cite{Frankel2}, and \cite{Galloway} to our attention.

\section{Theorem \ref{theorem1.3} and consequences}
\setcounter{equation}{0}

\subsection{The proof}
The proof of Theorem \ref{theorem1.3} is by a standard application of the second variation formula.
See \cite{Frankel2, PW} for the $m=0$ case and \cite[Theorem 7.4]{WW} for the $m=\infty$, $X=\nabla f$ case.

\begin{proof}[Proof of Theorem \ref{theorem1.3}]
Let $\alpha:[-\varepsilon,\varepsilon]\times [0,\ell]$ be a variation of a unit speed curve $\gamma(\cdot):=\alpha (0,\cdot)$ and define the arclength of the curve $\gamma_s(t):=\alpha(s,t)$ by
\begin{align}
\label{eq2.1}
L(s)=\int\limits_0^\ell\left \vert \frac{\partial \alpha}{\partial t}(s,t)\right \vert dt\ .
\end{align}
By Synge's second variation formula, whenever $\gamma(\cdot):=\alpha(0,\cdot)$ is a stationary point of $L(s)$ then
\begin{equation}
\label{eq2.2}
\frac{d^2L}{ds^2}(0)=\int\limits_0^\ell \left [ \left \vert \nabla_{\gamma'} V\right \vert^2 -\scc(V,\gamma')\right ] dt +\left \langle \gamma',\brax{\nabla_{V}\pdx[\alpha]{s}}(0,t)\right \rangle \bigg\vert_0^\ell\ ,
\end{equation}
where $V(t)=\frac{\partial \alpha}{\partial s}(0,t)$ is the variation vector field, with $s$ chosen here so that $V$ is orthogonal to $\gamma'(0)$ along $\gamma$.

Let $N_1$, $N_2\subset M$ be hypersurfaces and let $p_i\in N_i$ be the points in the hypersurfaces which are closest to each other. By way of contradiction, assume that $p_1\neq p_2$ so that $\dist(N_1,N_2)=:\ell>0$. Then choose a unit speed minimizing geodesic $\gamma:[0,\ell]\to M$ from $p_1$ to $p_2$. Next, select an orthonormal frame at $p_1$ such that the $n^{\rm th}$ element is $\gamma'(0)$ and parallel-transport it along $\gamma$ to define a basis $\{ E_1,...\, ,E_{n-1}, E_n\}$ along $\gamma$ with $E_n=\gamma'$. At the endpoints of $\gamma$, the basis vectors $E_1,...\, ,E_{n-1}$ will be tangent to the hypersurfaces. Pick variations $\alpha_1,...\, ,\alpha_{n-1}$ with the property that $\alpha_j(s,0)\in N_1$, $\alpha_j(s,\ell)\in N_2$ for sufficiently small $s$, and such that
\begin{align}
\label{eq2.3}
\pdx[\alpha_j]{s}(0,t)=E_j\ .
\end{align}
Then $\nabla_{\gamma'} V  = \nabla_{\gamma'} E_j =0$, so the first term on the right of \eqref{eq2.2} vanishes. Also, let $L_j$ denote the length on the left-hand side of \eqref{eq2.1} when $\alpha_j$ replaces $\alpha$ on the right. Summing the contributions of all $n-1$ such variations, then
\begin{equation}
\label{eq2.4}
\begin{split}
\sum\limits_{j=1}^{n-1}\frac{d^2L_j(0)}{ds^2}
=&\, -\sum\limits_{j=1}^{n-1}\int\limits_0^\ell \scc (E_j,\gamma') dt +\sum\limits_{j=1}^{n-1}\left \langle \gamma', \brax{\nabla_{\pdx[\alpha_j]{s}} \pdx[\alpha_j]{s}}(0,t)\right \rangle \at_0^\ell\\
=&\,-\int\limits_0^\ell\ric(\gamma',\gamma')dt +H(\ell)+H(0)
\end{split}
\end{equation}
by Synge's formula, where $H(0)$ denotes the mean curvature of $N_1$ at $p_1$ with respect to $\gamma'(0)$ and $H(\ell)$ denotes the mean curvature of $N_2$ at $p_2$ with respect to $-\gamma'(\ell)$, this orientation being consistent with our conventions in the case when $N_1$ and $N_2$ are boundary hypersurfaces for an interior region containing $\gamma$. We have
\begin{equation}
\label{eq2.5}
\begin{split}
\sum\limits_{j=1}^{n-1}\frac{d^2L_j(0)}{ds^2}
=&\, -\int\limits_0^\ell\ric_X^m(\gamma',\gamma')dt+
\int\limits_0^\ell\left [ \frac{1}{2}\pounds_Xg(\gamma',\gamma') -\frac{1}{m}\left ( g(X,\gamma')\right )^2\right ] dt
+ H(\ell) + H(0)\\
= &\,  -\int\limits_0^\ell\ric_X^m(\gamma',\gamma')dt + g(X,\gamma')\at_{0}^{\ell} -\frac{1}{m}\int\limits_0^\ell \left ( g(X,\gamma')\right )^2 dt  +H(\ell) +H(0) \\
=&\,  -\int\limits_0^\ell\ric_X^m(\gamma',\gamma')dt +H_X(\ell) +H_X(0)  -\frac{1}{m}\int\limits_0^\ell \left ( g(X,\gamma')\right )^2 dt\\
<&\, -\frac{1}{m}\int\limits_0^\ell \left ( g(X,\gamma')\right )^2 dt\ .
\end{split}
\end{equation}
where in the last step we used our curvature assumptions. All steps are valid for any $m$, including negaitve $m$ and $m=\infty$ provided $1/m$ is interpreted as $0$ in that case.

Thus, whenever $m\in (0,\infty]$ the base geodesic $\gamma(t)=\alpha(0,t)$ must be unstable. This contradicts the assumption that $\gamma$ is a minimizing curve between closed hypersurfaces.
\end{proof}

The proof shows that stronger results can be obtained. For example, the last step of \eqref{eq2.5} follows provided only that $\int\limits_0^\ell\ric_X^m(\gamma',\gamma')dt -H_X(\ell)- H_X(0)>0$ along each unit speed geodesic $\gamma$. For finite $m$, we may replace the open inequality ($<0$) by the closed one ($\le 0$) in \eqref{eq2.5} unless $X$ and $\gamma$ are orthogonal all along $\gamma$.

\subsection{Corollaries and extensions}
If $\ric_X^m\ge kg$ for some $k>0$, then the fundamental group of $M$ is finite, just as in the $m=0$ case. But an infinite fundamental group may arise when $\ric_X^m>0$ if there is no such $k$ and if $M$ is not closed, when there is also no positive lower bound on the induced Bakry-\'Emery Ricci curvature of the $X$-minimal surface $N$. But in general, theorem \ref{theorem1.3} has the following topological consequence when $\ric_X^m>0$, the $m=0$ case of which was proved in \cite{Frankel2}.

\begin{corollary}\label{corollary2.1}
Let $M$ be a complete manifold with a vector field $X$ and an
$m\in (0,\infty]$ such that $\ric_X^m>0$ pointwise on $M$, and such that $M$ has a closed Bakry-\'Emery $X$-minimal surface $N$. Then the homomorphism $i_*:\pi_1(N)\to\pi_i(M)$ induced by inclusion is surjective.
\end{corollary}

\begin{proof} Suppose not. Then there is a class $[c']\in\pi_1(M,p)$ of loops based at $p\in N$ which cannot be deformed to lie in $N$, and for which $c'$ is a length-minimizing representative. If we further minimize over all $p$ in the closed submanifold $N$, there will be a shortest nontrivial loop $c$ of length $L(c)>0$. We may pass to the universal covering space ${\bar M}$, with Riemannian metric ${\bar g}$ and vector field ${\bar X}$ on ${\bar M}$ defined by pullback. Then the $m$-Bakry-\'Emery Ricci tensor on ${\bar M}$ will obey $\ric_{\bar X}^m({\bar g})>0$, and there will be two disjoint copies of the Bakry-\'Emery $X$-minimal surface $N$ joined by a minimal geodesic of length $L(c)>0$. But by Theorem \ref{theorem1.3} these two hypersurfaces must intersect, which is a contradiction.
\end{proof}

The above corollary would also hold if the pointwise assumption $\ric_X^m>0$ were replaced by the condition that $\int\limits_0^\ell\ric_X^m(\gamma',\gamma')dt>0$ on each closed geodesic loop $c$.

Finally, we give an extension of Theorem \ref{theorem1.3} to the case where the Bakry-\'Emery Ricci curvature bound may be negative. We state and prove this result only for finite $m$.

\begin{lemma}\label{lemma2.3}
Let $M$ be a compact manifold with vector field $X$ and an $m\in (0,\infty)$ such that for some $k> 0$ we have $\ric_X^m > -(n+m-1)k$ pointwise on $M$. Let $M$ contain a hypersurface $N$ that bounds a connected region $\Omega$ of compact closure ${\bar \Omega}$ in $M$.
Suppose that the Bakry-\'Emery mean curvature of $N$ (defined by the normal pointing into $\Omega$) satisfies $H_X\le -(n+m-1)\sqrt{k}$. Then $N$ is connected, and the map $i_*:\pi_1(N)\to\pi_1(\bar \Omega)$ induced by inclusion is surjective.
\end{lemma}

\begin{proof}
Suppose by way of contradiction that $N_1$ and $N_2$ are distinct non-empty connected components of $N$. Once again find the unit speed geodesic that realizes the distance between these components. We let $\phi$ be the solution of $\phi''-k\phi=0$, with $\phi(0)=\phi(\ell)=1$; i.e., $\phi(x) =\frac{\cosh(\sqrt{k} (x-\ell/2))}{\cosh (\sqrt{k}\ell/2)}$. In \eqref{eq2.1} we let $V(t)=\phi(t)E_i(t)$. Then \eqref{eq2.2} becomes
\begin{equation}
\label{eq2.8}
\frac{d^2L_j(0)}{ds^2} = \int\limits_0^{\ell} \left [ \phi'^2-\phi^2\scc(E_j,\gamma')\right ] dt +\phi^2(0,t)\left \langle \gamma',\nabla_{E_i}E_i\right \rangle \bigg\vert_0^\ell\ .
\end{equation}
We sum over $i\in \{ 1,\dots,n-1\}$ to get
\begin{equation}
\label{eq2.9}
\begin{split}
\sum\limits_{j=1}^{n-1}\frac{d^2L_j(0)}{ds^2}
=&\, \int\limits_0^\ell \left [ (n-1)(\phi')^2 -\phi^2\ric(\gamma',\gamma')\right ] dt +\sum_{i=1}^{n-1}\phi^2\metr{\nabla_{E_i}E_i,\gamma'}\at_0^\ell\\
=&\, \int\limits_0^\ell \left [ (n-1)(\phi')^2 -\phi^2\ric_X^m(\gamma',\gamma') + \frac{\phi^2}{2} \pounds_X g(\gamma',\gamma')-\frac{\phi^2}{m}\left \langle X,\gamma'\right \rangle^2\right ] dt \\
&\, +H(\ell)+H(0).
\end{split}
\end{equation}
We simplify the first two terms on the right as follows.
\begin{equation}
\label{eq2.10}
\begin{split}
\int\limits_0^\ell \left [ (n-1)(\phi')^2 -\phi^2\ric_X^m(\gamma',\gamma')\right ] dt = &\, (n-1)\phi\phi'\big\vert_0^{\ell} -\int\limits_0^{\ell}\left [ (n-1)\phi\phi'' +\phi^2\ric_X^m(\gamma',\gamma')\right ] dt\\
=&\, 2(n-1)\sqrt{k}\tanh\frac{\sqrt{k}\ell}{2} -\int\limits_0^{\ell}\phi^2\left [ (n-1)k +\ric_X^m(\gamma',\gamma')\right ] dt\\
< &\, 2(n-1)\sqrt{k}\tanh\frac{\sqrt{k}\ell}{2}+mk\int\limits_0^{\ell}\phi^2dt.
\end{split}
\end{equation}
Furthermore, the third term on the right of \eqref{eq2.9} is
\begin{equation}
\label{eq2.11}
\int\limits_0^\ell \frac{1}{2}\phi^2 \pounds_Xg(\gamma',\gamma')dt = \int\limits_0^\ell \phi^2 \nabla_{\gamma'}\left \langle X,\gamma'\right \rangle dt = \left \langle X,\gamma'(\ell)\right \rangle -\left \langle X,\gamma'(0)\right \rangle - 2\int\limits_0^\ell \phi\phi' \left \langle X,\gamma'\right \rangle dt,
\end{equation}
where we used that $\phi(\ell)=\phi(0)=1$. Since we define $H$ using a normal vector parallel to $\gamma'$ at $t=0$ and anti-parallel to it at $t=\ell$, then $H_X(\ell)=H(\ell)+\left \langle X,\gamma'(\ell)\right \rangle$ and $H_X(0)=H(0)-\left \langle X,\gamma'(0)\right \rangle$, so we may combine the $H$-terms in \eqref{eq2.9} with the $\left \langle X,\gamma'\right \rangle$ terms in \eqref{eq2.11}. Putting all this together, \eqref{eq2.9} becomes
\begin{equation}
\label{eq2.12}
\begin{split}
\sum\limits_{j=1}^{n-1}\frac{d^2L_j(0)}{ds^2}< &\, 2(n-1)\sqrt{k}\tanh\frac{\sqrt{k}\ell}{2} +mk\int\limits_0^{\ell}\phi^2dt -\int\limits_0^\ell \left [ 2 \phi\phi' \left \langle X,\gamma'\right \rangle +\frac{1}{m}\phi^2  \left ( g(X,\gamma')\right )^2\right ] dt\\
&\, + H_X(\ell) + H_X(0)\\
=&\, 2(n-1)\sqrt{k}\tanh\frac{\sqrt{k}\ell}{2} +m\int\limits_0^{\ell}\left [ k\phi^2+\phi'^2\right ]dt -\int\limits_0^\ell \left [ \sqrt{m}\phi' +\frac{\phi}{\sqrt{m}}\left ( g(X,\gamma')\right )\right ]^2 dt\\
&\, + H_X(\ell) + H_X(0)\\
=&\, -\int\limits_0^{\ell} \left [ \sqrt{m}\phi'+\frac{\phi}{\sqrt{m}}g(X,\gamma')\right ]^2 dt+ H_X(\ell) + H_X(0)+ 2(n+m-1)\sqrt{k}\tanh\frac{\sqrt{k}\ell}{2}.
\end{split}
\end{equation}
Using $H_X(t)\le -(n+m-1)\sqrt{k}$ for $t=0$, $\ell$, then
$\sum\limits_{j=1}^{n-1}\frac{d^2L_j(0)}{ds^2} <0$,
which is a contradiction, so $N$ must be connected. Similarly, $\pi^{-1}(N)$ is connected in the universal cover $\pi:{\tilde M}\to M$ of $M$. Let $\ell$ be a loop in $\bar \Omega$ based at $p\in N$. Lift this to a curve ${\tilde  \ell}$ in ${\tilde {\bar \Omega}}\subset{\tilde M}$ joining points $p_1$, $p_2\in \pi^{-1}(p)$. Let ${\tilde \ell}'$ be a curve joining $p_1$ to $p_2$, contained in $\pi^{-1}(N)$. Since ${\tilde M}$ is simply connected, ${\tilde \ell}$ and ${\tilde \ell}'$ are homotopic, so $\pi({\tilde \ell})=\ell$ and $\pi({\tilde \ell}')$ are as well. Thus $i_*$ is surjective.
\end{proof}

\section{$X$-mean curvature of level sets of distance functions}
\setcounter{equation}{0}

\noindent
\subsection{The Bochner formula} We define the weighted Laplacian on functions to be
\begin{align}
\label{eq3.1}
\Delta_X \phi = \Delta \phi - g\brax{X,\nabla \phi}\ ,
\end{align}
where $\Delta\phi :=\tr_g \hess \phi$. We also define $\tf \hess d:=\hess d -\frac{1}{(n-1)}h\Delta d$ to be the $h$-tracefree part of $\hess d$, where $h$ is the induced metric on level sets of $d$.

The Bakry-\'Emery version of the Bochner identity (e.g., \cite[Lemma 4]{KWW} which is concerned with finite positive $m$, but the identity will hold for any $m\neq 0$, including $m=\infty$ if in that case $1/m$ is interpreted as $0$), is given by
\begin{align}
    \label{eq3.2}
    \Delta_X(\abs{\nabla w}^2)=2\abs{\hess w}^2 + 2\nabla_{\nabla w}(\Delta_X w)+2\ric_X^m(\nabla w,\nabla w)+\frac{2}{m}[X(w)]^2\ .
\end{align}
Let $d$ be a distance function, so $|d|=1$. Application of \eqref{eq3.2} with $w=d$ yields
\begin{equation}
\label{eq3.3}
\begin{split}
\nabla_{\nabla d}\Delta_X d=&\, -\left \vert \hess d \right \vert^2 - \ric_X^m(\nabla d, \nabla d) - \frac{1}{m}[X(d)]^2\\
=&\,  -\left \vert \tf \hess d+\frac{1}{(n-1)}\left [ \Delta_X d +(X(d))^2\right ] \right \vert^2 -\ric_X^m(\nabla d, \nabla d) - \frac{1}{m}[X(d)]^2\\
=&\, -\left \vert \tf \hess d\right \vert^2 - \ric_X^m(\nabla d, \nabla d) -\frac{\left ( \Delta_X d\right )^2}{(n+m-1)}\\
&\, -\frac{1}{(n-1)}\left ( \sqrt {\frac{m}{n+m-1}}\Delta_X d +\sqrt{\frac{n+m-1}{m}}X(d)\right )^2 \ .
\end{split}
\end{equation}

\subsection{The Riccati equation} Let $\gamma$ belong to a congruence unit-speed geodesics parametrized by $t$, either issuing from a point or issuing orthogonally from an initial hypersurface $\Sigma$ at $t=0$. Since such a congruence is irrotational, we can consider the hypersurfaces defined by level sets $t=d$ at equal parameter values along different curves in the congruence. The mean curvature $H(t)$ of these level sets is governed by the scalar Riccati equation (sometimes called the Raychaudhuri equation)
\begin{equation}
\label{eq3.4}
\frac{dH}{dt}=-\ric(\gamma',\gamma')-|A|^2=-\ric(\gamma',\gamma')-|\tf A|^2-\frac{H^2}{(n-1)} \ ,
\end{equation}
where $A:=\hess w$ denotes the second fundamental form of the hypersurface and $H$ is its trace. This equation can be rewritten in terms of Bakry-\'Emery quantities as
\begin{equation}
\label{eq3.5}
\begin{split}
&\, \ddt[H]-\frac{1}{2}\calL_Xg(\gamma',\gamma') -\ric_X^m(\gamma',\gamma')-|\tf A|^2 -\frac{1}{m}\left ( g(\gamma',X) \right )^2 -\frac{H^2}{(n-1)}\\
\implies &\, \frac{dH_X}{dt}= -\ric_X^m(\gamma',\gamma')-|\tf A|^2-\frac{1}{m}\left ( g(\gamma',X)\right )^2 -\frac{2}{(n-1)}H_Xg(\gamma',X)-\frac{H_X^2}{(n-1)}-\frac{\left ( g(\gamma',X)\right )^2}{(n-1)}\\
\implies &\, \frac{dH_X}{dt}= -\ric_X^m(\gamma',\gamma')-|\tf A|^2-\frac{H_X^2}{(n+m-1)} -\frac{1}{(n-1)}\left [ \sqrt{\frac{m}{n+m-1}} H_X+\sqrt{\frac{n+m-1}{m}}g(\gamma',X)\right ]^2
\end{split}
\end{equation}
which is just equation \eqref{eq3.3}. As before, the $m\to\infty$ of this expression will yield the correct $m=\infty$ result.

\subsection{Consequences I: Distance from a hypersurface} Using the above equations, we are able to give bounds on the $X$-mean curvature of level sets of distance functions that realize the distance from a given hypersurface.

\begin{lemma}[Finite $m>0$] \label{lemma3.1}
Let there be an $X$, an $m\in (0,\infty)$, and a $\delta\in \{ 0,1\}$ such that $\ric_X^m\ge -(n-1)\delta g$. Let $H_X(t)$ denote the Bakry-\'Emery mean curvature of the $t$-level set $d(\Sigma,\cdot)=t$ of the distance function from the initial hypersurface $\Sigma$ (at $t=0$), as long as this function is smooth. If
\begin{itemize}
\item[(i)] $\delta=0$, or
\item[(ii)] $\delta=1$ and $(H_X(0))^2\ge (n-1)(n+m-1)$,
\end{itemize}
then $H_X(t)\le H_X(0)$ for all $t>0$ for which $H_X$ is defined. Furthermore, we have the following.
\begin{itemize}
\item[a)] If $H_X(t_1)=H_X(0)\le 0$ for any $t_1>0$ then $H_X(t)=H_X(0)$ for all $0\le t\le t_1$ and then along a minimizing geodesic $\gamma:[0,t_1]\to M$ we have that $H_X(0)=-\sqrt{(n-1)(n+m-1)}\delta=H_X(t)$, $\ric_X^m(\gamma',\gamma')=-(n-1)\delta g$, the tracefree part $\tf A$ of the second fundamental form vanishes, and $g(\gamma',X)=- m\delta \sqrt{\frac{n-1}{n+m-1}}$.
\item[b)] If $H_X(0)<-\delta \sqrt{(n-1)(n+m-1)}$ then the domain of $\gamma$ is bounded.
\end{itemize}
\end{lemma}

\begin{proof}
Defining $x(t):=\frac{H_X}{\sqrt{(n-1)(n+m-1)}}$ and $c:=\sqrt{\frac{n-1}{n+m-1}}$, then equation \eqref{eq3.5} yields the inequality
\begin{equation}
\label{eq3.6}
x'(t)\le c\left ( \delta -x^2(t)\right )\ .
\end{equation}
Let $u(t):= e^{c\int\limits_0^t x(s) ds}$. Then $x= \frac{u'}{cu}$, and the inequality \eqref{eq3.6} becomes $\frac{u''}{cu}\le c\delta $, or $u''\le c^2\delta u$.

Let $v(t)$ be the unique solution of $v''=c^2\delta v$ such that $v(0)=u(0)=:u_0>0$ and $v'(0)=u'(0)=:u_0'$. Let $T>0$ be the first point at which either $u(t)$ or $v(t)$ has a zero or becomes undefined, and restrict attention to $t\in [0,T)$. Now let $y(t):=\frac{v'}{cv}$. A simple calculation shows that $[uv(x-y)]' =\frac{1}{c}(vu''-uv'') \le 0$. Integrating and using $x(0)=y(0)$, we have that $uv(x-y)\le 0$.  Then $x(t)\le y(t)$ for all $t\in[0,T)$.

When $\delta=0$, then $v(t)=u_0't+u_0$ and so
\begin{equation}
\label{eq3.7}
x(t)\le y(t)=\frac{u_0'}{c(u_0't+u_0)} =\frac{u_0'/(cu_0)}{\frac{u'_0}{u_0}t+1} =\frac{x(0)}{cx(0)t+1}\le x(0)\text{\ for\ }t\in [0,T).
\end{equation}
When $\delta=1$, then $v(t)=u_0\cosh ct +\frac{u_0'}{c}\sinh ct$. If $[x(0)]^2\ge 1$ as well, then we have
\begin{equation}
\label{eq3.8}
x(t)\le y(t)=\frac{u_0\sinh ct +\frac{u_0'}{c}\cosh ct}{u_0\cosh ct +\frac{u_0'}{c}\sinh ct} =\frac{\tanh ct +\frac{u_0'}{cu_0}}{1 +\frac{u_0'}{cu_0}\tanh ct} =\frac{x(0)+\tanh ct }{1 +x(0)\tanh ct}\le x(0)\text{\ for\ }t\in [0,T).
\end{equation}
This proves the inequality $H_X(t)\le H_X(0)$ for all $t>0$. It also proves part (b), since for the $\delta=0$ case if $x(0)<0$ in \eqref{eq3.7} then $x(t)\to -\infty$ at or before $t_1:=-\frac{1}{cx(0)}>0$, while for $\delta=1$ the corresponding result follows from taking $x(0)<-1$ in \eqref{eq3.8}.

To prove the equality statement, observe that the last (i.e., rightmost) inequalities in \eqref{eq3.7} and \eqref{eq3.8} are strict for $t>0$ unless $x(0)=0$ in \eqref{eq3.7} or $x(0)=\pm 1$ in \eqref{eq3.8}, and then $H_X(0)=\pm \sqrt{(n-1)(n+m-1)}\delta$. So we need only consider the case of $H_X(t)=- \sqrt{(n-1)(n+m-1)}\delta$ for some $t>0$. But then there is a local minimum of $H_X$ at some $0<\tau <t$ such that $H_X'(\tau)=0$. Since the left-hand side of \eqref{eq3.5} vanishes there, the right-hand side must as well. But under the given conditions we have $-\ric_X^m(\gamma',\gamma')-\frac{H_X^2}{(n+m-1)}\le 0$. Then $\ric_X^m(\gamma',\gamma')=-\frac{H_X^2}{(n+m-1)}=-(n-1)\delta$ at $\tau$, and so $H_X(t)=-\sqrt{(n-1)(n+m-1)}\delta$ (the minimum condition eliminates the positive root) at $\tau$ and so for all $t$. But then $H_X(t)$ is constant and $dH/dt=0$ for all $t$, so $\ric_X^m(\gamma',\gamma')+\frac{H_X^2}{(n+m-1)}$, $|\tf A|$, and the final term on the right in square brackets must each vanish independently throughout the domain. This can only happen when the conditions listed in the theorem hold.
\end{proof}

Furthermore, this result extends to $m=\infty$ and $m<1-n$ when $\delta=0$.
\begin{lemma}[Case of $m=\infty$]\label{lemma3.2}
Let there be an $X$ such that $\ric_X^{\infty}\ge 0$. Let $H_X(t)$ denote the Bakry-\'Emery mean curvature of the $t$-level set $d(\Sigma,\cdot)=t$ of the distance function from the initial hypersurface $\Sigma$ (at $t=0$), as long as this function is smooth. Then $H_X(t)\le H_X(0)$ for all $t>0$ for which $H_X$ is defined. If $H_X(t_1)=H_X(0)\le 0$ for any $t_1>0$ in the domain then $H_X(t)=H_X(0)=-g(\gamma',X)$ for all $0\le t\le t_1$ along a minimizing geodesic $\gamma:[0,t_1]\to M$, and $\ric_X^m(\gamma',\gamma')=0$. The tracefree part $\tf A$ of the second fundamental form also vanishes.

Finally, if $\int\limits_0^t g(\gamma'(\tau),X)d\tau \le K$ for some constant $K>0$ and if $H_X(0)<0$, then the domain of $\gamma$ is bounded.
\end{lemma}

\begin{proof}
In this case, equation \eqref{eq3.5} yields
\begin{equation}
\label{eq3.9}
\frac{dH_X}{dt}=-\ric_X^{\infty}(\gamma',\gamma')-\left \vert \tf A\right \vert^2-\frac{1}{(n-1)}\left ( H_X+g(\gamma',X)\right )^2.
\end{equation}
Since $\ric_X^{\infty}\ge 0$, each term on the right is negative semi-definite, so $H_X(t)\le H_X(0)$ for $t>0$, and $H_X(t_1)=H_X(0)$ if and only if $\ric_X^{\infty}(\gamma',\gamma')=0$, ${\tf A}=0$, and $H_X+g(\gamma',X)=0$ for all $0\le t\le t_1$. But then $H_X(0)=H_X(t)=-g(\gamma',X)(t)=-g(\gamma',X)(0)$ for all $0\le t\le t_1$.

It remains to prove the boundedness of the domain when $H_X(0)<0$. The $m=\infty$ case of equation \eqref{eq3.5} can be written as
\begin{equation}
\label{eq3.10}
\begin{split}
\frac{dH_X}{dt}+\frac{2}{(n-1)} g(\gamma',X) H_X = &\, -\ric_X^m(\gamma',\gamma') -|\tf A|^2 -\frac{H_X^2}{(n-1)} -\frac{n+m-1}{m(n-1)}\left (g(\gamma',X)\right )^2\\
\le &\, -\frac{H_X^2}{(n-1)}.
\end{split}
\end{equation}
Along $\gamma$, define
\begin{equation}
\label{eq3.11}
{\mathcal H}_X(t):= e^{\frac{2}{(n-1)}\int\limits_0^t g(\gamma', X)ds}H_X(t),
\end{equation}
then \eqref{eq3.10} becomes
\begin{equation}
\label{eq3.12}
\frac{d{\mathcal H}_X}{dt}\le  - \frac{e^{-\frac{2}{(n-1)}\int\limits_0^t g(\gamma', X)ds}}{(n-1)}{\mathcal H}_X^2 \le \frac{e^{-\frac{2K}{(n-1)}}}{(n-1)}{\mathcal H}_X^2 .
\end{equation}
Since ${\mathcal H}_X(0)<0$, there is an interval beginning at $t=0$ for which ${\mathcal H}_X(t)\neq 0$. On the interval, we can divide \eqref{eq3.12} by $\left ( {\mathcal H}_X(t)\right )^2$ and integrate to obtain
\begin{equation}
\label{eq3.13}
\frac{1}{{\mathcal H}_X(t)}\ge \frac{1}{{\mathcal H}_X(0)}+ \frac{e^{-\frac{2K}{(n-1)}}}{(n-1)}t.
\end{equation}
But then there will be some $t_1>0$ such that the right-hand side becomes zero, and therefore ${\mathcal H}_X(t)\to -\infty$ at some $t_2\in (0,t_1]$.
\end{proof}

\begin{lemma}[Case of $m\le 1-n$]\label{lemma3.3}
For $m\in (-\infty, 1-n]$, let there be an $X$ such that $\ric_X^m\ge 0$. Let $H_X(t)$ denote the Bakry-\'Emery mean curvature of the $t$-level set $d(\Sigma,\cdot)=t$ of the distance function from the initial hypersurface $\Sigma$ (at $t=0$), as long as this function is smooth, and let $\gamma$ be a minimizing geodesic. Then $e^{\frac{2}{n-1}\int_0^t g(\gamma',X)dr}H_X(t)\le H_X(0)$ for all $t>0$ for which $H_X$ is defined. If $e^{\frac{2}{n-1}\int_0^{t_1} g(\gamma',X)dr}H_X(t_1)=H_X(0)\le 0$ for some $t_1>0$ then $H_X(t)=H_X(0)=0$ for all $0\le t\le t_1$ and either $g(\gamma',X)=0$ for all $0\le t\le t_1$ as well or $m=1-n$. Furthermore, along $\gamma:[0,t_1]\to M$ we have that $\ric_X^m(\gamma',\gamma')=0$ and $\tf A=0$ .

Furthermore, if $\int\limits_0^t g(\gamma'(\tau),X)d\tau \le K$ for some $K>0$ and if $H_X(0)<0$ then the domain of $\gamma$ is bounded.
\end{lemma}

\begin{proof} Equations \eqref{eq3.10}--\eqref{eq3.13} hold when $m\le 1-n$, proving the claim about the boundedness of the domain.

Next, the middle equality in \eqref{eq3.5} yields
\begin{equation}
\label{eq3.14}
\frac{dH_X}{dt}= -\ric_X^m(\gamma',\gamma')-\abs{\tf A}^2-\frac{1}{m}g(\gamma',X)^2-\frac{2}{n-1}H_Xg(\gamma',X) -\frac{H_X^2}{(n-1)}-\frac{1}{(n-1)}\left ( g(\gamma',X) \right )^2,
\end{equation}
so
\begin{equation}
\label{eq3.15}
e^{\frac{-2}{n-1}\int g(\gamma',X)dr}\brax{e^{\frac{2}{n-1}\int g(\gamma',X)dr}H_X}' =
-\ric_X^m(\gamma',\gamma')-\abs{\tf A}^2-\frac{H_X^2}{(n-1)}-\frac{(n+m-1)}{m(n-1)}\left ( g(\gamma',X) \right )^2.
\end{equation}
Since $\ric_X^{m}\ge 0$, then each term on the right is negative semi-definite so we have
\begin{equation}
\label{eq3.16}
e^{\frac{2}{n-1}\int_0^s g(\gamma',X)dr}H_X(s)\leq e^{\frac{2}{n-1}\int_0^t g(\gamma',X)dr}H_X(t)\le H_X(0)
\end{equation}
for $s> t>0$, and $e^{\frac{2}{n-1}\int_0^{t_1} g(\gamma',X)dt}H_X(t_1)=H_X(0)$ if and only if $\ric_X^{m}(\gamma',\gamma')=0$, ${\tf A}=0$, $g(\gamma',X)=0$, and $H_X(t)=0$ for all $0\le t\le t_1$. Then $H_X(0)=H_X(t)=0$ for all $0\le t\le t_1$, and either $g(\gamma',X)=0$ for all $0\le t\le t_1$ as well or $m=1-n$.
\end{proof}

\begin{remark}\label{remark3.4}
The prefactors $e^{\pm\frac{2}{n-1}\int_0^t g(\gamma',X)dr}$ in various expressions above are path-dependent in general but are well-defined here, since the path $\gamma$ is uniquely specified.
\end{remark}

\begin{remark}\label{remark3.5} Lemmata \ref{lemma3.2} and \ref{lemma3.3} generalize only the $\delta=0$ case of Lemma \ref{lemma3.1}. We were not able to generalize the $\delta=1$ case for $m=\infty$ or $m\in(-\infty,1-n]$.
\end{remark}

\subsection{Consequences II: Distance from a point} We now modify the arguments of the previous subsection to obtain information about the mean curvature of level sets of distance functions that realize the distance from a point.

\begin{lemma}[Finite $m>0$] \label{lemma3.6}
Let there be an $X$, an $m\in (0,\infty)$, and a $\delta\in \{ 0,1\}$ such that $\ric_X^m\ge -(n-1)\delta g$. Let $H_X(t)$ denote the Bakry-\'Emery mean curvature of the $t$-level set $d(p,\cdot)=t$ of the distance function from a chosen initial point $p$ (at $t=0$), as long as this function is smooth. Let $v(t)$ denote the unique solution of $v''(t)-\delta v(t)=0$ with $v(0)=0$ and $v'(0)=1$. Let $c:=\sqrt{\frac{n-1}{n+m-1}}$ and define $y(t):=\frac{v'(t)}{cv(t)}$ for $t>0$. Then
\begin{equation}
\label{eq3.17}
H_X(t)\le \sqrt{(n-1)(n+m-1)}y(t)=(n-1)\begin{cases} 1/t,& \delta=0,\\ \coth{t}, & \delta=1, \end{cases}
\end{equation}
for all $t>0$ for which $H_X$ is defined. Equality holds iff $\ric_X^m(\gamma',\gamma')= -(n-1)\delta$, $\tf A=0$, and $g(\gamma',X)=-\frac{mH_X}{(n+m-1)}=-\frac{m(n-1)}{(n+m-1)}\begin{cases} 1/t,& \delta=0,\\ \coth{t}, & \delta=1, \end{cases}$ for all $t>0$ in the domain.
\end{lemma}

From the proof of Lemma \ref{lemma3.1}, we see that under the conditions of Lemma \ref{lemma3.6} then $v(t) = t$ if $\delta=0$ and $v(t)=\frac{1}{c}\sinh (ct)$ when $\delta=1$.

\begin{proof}
Choose $0<t_0\le t$. Define $x(t):=\frac{H_X}{\sqrt{(n-1)(n+m-1)}}$ and $u(t):= e^{c\int\limits_{t_0}^t x(s) ds}$. Then the simple calculation that appears in the proof of Lemma \ref{lemma3.1} again yields that $[uv(x-y)]' \le 0$, though only for $t\ge t_0>0$ since $y(t)$ is now not defined at $t=0$. We integrate on the domain $0<t_0\le t$ to obtain $u(t)v(t)(x(t)-y(t))\le u(t_0)v(t_0)(x(t_0)-y(t_0)) = \frac{1}{c}\left ( v(t_0)u'(t_0)-u(t_0)v'(t_0)\right )$. Now take $t_0\searrow 0$ and use that $u(0)=v(0)=0$ to obtain, as in Lemma \ref{lemma3.1}, that $u(t)v(t)(x(t)-y(t))\le 0$. In either case, then $x(t)\le y(t)$. Then \eqref{eq3.17} follows.

For the equality statement, to prove the forward implication (``if'') we plug $\ric_X^m(\gamma',\gamma')= -(n-1)\delta$, $\tf A=0$, and $g(\gamma',X)=-\frac{mH_X}{(n+m-1)}$ into the last line of \eqref{eq3.5} to obtain the differential equation $H_X'(t)=(n-1)\delta -H_X^2/(n+m-1)$, which is easily solved to find $H_X$ as claimed in the lemma.

For the reverse implication (``only if''), substitute $H_X(t)=(n-1)\begin{cases} 1/t,& \delta=0,\\ \coth{t}, & \delta=1, \end{cases}$ in the last line of \eqref{eq3.5}.
\end{proof}

\begin{lemma}[$m=\infty$ and $m\le 1-n$] \label{lemma3.7}
Let there be an $X$ such that $\ric_X^{\infty}\ge 0$ and an $m\in \{\infty\}\cup (-\infty,1-n]$ such that $\left \vert \int_{t_0}^m g(\gamma', X)dt\right \vert\le K$ along every minimizing geodesic $\gamma:[t_0,\infty)\to M$. Let $H_X(t)$ denote the Bakry-\'Emery mean curvature of the $t$-level set $d(p,\cdot)=t$ of the distance function from a chosen initial point $p$ (at $t=0$), as long as this function is smooth. Then for all $t>0$ we have
\begin{equation}
\label{eq3.18}
H_X(t)\le \frac{(n-1)e^{4K/(n-1)}}{t}.
\end{equation}

\end{lemma}

\begin{proof}
If $m=\infty$, then from \eqref{eq3.9}, we can write
\begin{equation}
\label{eq3.19}
\frac{d}{dt} \left [ e^{\frac{2}{(n-1)}\int\limits_0^t g(\gamma', X)ds }H_X(t) \right ] =-e^{\frac{2}{(n-1)}\int\limits_0^t g(\gamma', X)ds } \left \{ \ric_X^{\infty} +|\tf A|^2 + \frac{\left [ H_X^2 +\left ( g(\gamma',X)\right )^2\right ])}{(n-1)}\right \}.
\end{equation}
Let ${\mathcal H}_X(t):= e^{\frac{2}{(n-1)}\int\limits_0^t g(\gamma', X)ds}H_X(t)$. Then we have
\begin{equation}
\label{eq3.20}
\frac{d}{dt} {\mathcal H}_X(t) \le -e^{-\frac{2}{(n-1)}\int\limits_0^t g(\gamma'(s), X)ds} {\mathcal H}_X^2(t).
\end{equation}
If $m\in(-\infty,1-n]$, then \eqref{eq3.15} also leads to \eqref{eq3.20}.

Integrating and simplifying, we eventually obtain
${\mathcal H}(t)\le \left ( \frac{1}{(n-1)}\int\limits_{t_0}^t e^{-\frac{2}{(n-1)}\int_0^{\tau} g(X,\gamma')ds}d\tau +\frac{1}{H(t_0)}\right )^{-1}$, so
\begin{equation}
\label{eq3.21}
H(t)\le (n-1)e^{-\frac{2}{(n-1)}\int_0^t g(X,\gamma'(r))dr}\left [ \int\limits_{t_0}^t e^{-\frac{2}{(n-1)}\int_0^{\tau} g(X,\gamma'(s))ds}d\tau +\frac{1}{H(t_0)}\right ]^{-1}.
\end{equation}
Now take $t_0\searrow 0$ and use that $H(t_0)\sim (n-1)/t_0 >0$ for small positive $t_0$ to get
\begin{equation}
\label{eq3.22}
H(t) \le (n-1)e^{-\frac{2}{(n-1)}\int_0^t g(X,\gamma'(r))dr}\left [\int\limits_0^t e^{-\frac{2}{(n-1)}\int_0^{\tau} g(X,\gamma'(s))ds}d\tau\right ]^{-1} .
\end{equation}
But the numerator is bounded above by $(n-1)e^{2K/(n-1)}$ and the denominator is bounded below by $\int_0^t (n-1)e^{-2K/(n-1)}dt = (n-1)e^{-2K/(n-1)}t$.
\end{proof}

\section{Topological results}
\setcounter{equation}{0}

\noindent
If the assumption in Theorem \ref{theorem1.3} is weakened to $\ric_X^m\ge 0$ and if there are now $X$-minimal surfaces $N_1$ and $N_2$ that do not intersect, there is a certain rigidity. For the $m=0$ case see \cite{Ichida}, \cite{Kasue}, \cite{Galloway}, \cite{CK}, and \cite[Theorem 4]{PW}. There are also topological implications when there is a single $X$-minimal surface if a locally isometric covering space can be found in which the single surface has multiple preimages under the covering map. For the $m=0$ case see, e.g., \cite[Theorem 2.5]{CF}. We repeat the argument, now taking $m\in (0,\infty]$. Reference \cite{CF} notes that some of the results follow from the Cheeger-Gromoll splitting theorem \cite{CG} if $m=0$; we note that for finite $m>0$ we can instead use the splitting theorem proved in \cite{KWW}. But we will give a proof directly from the Bochner formula.

\begin{lemma}\label{lemma4.1}
Let $M$ be a closed connected manifold with a vector field $X$ and an $m\in (0,\infty]$ such that $\ric_X^m \geq 0$ pointwise on $M$. Let $N_1$ and $N_2$ be two disjoint, closed, connected $X$-minimal hypersurfaces $N_1$ and $N_2$ such that $U$ is a connected component of $M\backslash \left ( N_1\cup N_2\right )$ with closure ${\bar U}\supset\left ( N_1\cup N_2\right )$. Let $d_i$ denote the distance map in $U$ from $N_i$, $i\in \{1,2 \}$.
Then $U$ splits isometrically as the product of an interval $I\subset {\mathbb R}$ with $N$, where $N$ is any level set of $d_i$ for either value of $i$ and all such level sets are isometric and totally geodesic. Furthermore, $\ric_X^m= 0$ and $X$ is tangent to level sets $N$ so $\Delta_X d_i=\Delta d_i=0$.
\end{lemma}

\begin{proof}
Equation \eqref{eq3.3} can be applied to each $d_i$ to yield
\begin{equation}
\label{eq4.1}
\begin{split}
0=&\, \left \vert \tf \hess d_i\right \vert^2 + \nabla_{\nabla d_i}(\Delta_X d_i)+\ric_X^m(\nabla d_i,\nabla d_i) +\frac{1}{(n-1)}\left ( \sqrt {\frac{m}{n+m-1}}\Delta_X d_i +\sqrt{\frac{n+m-1}{m}}X(d_i)\right )^2\ .
\end{split}
\end{equation}
The limit $m\to \infty$ of this expression gives the correct result for the $m=\infty$ case. In Lemma \ref{lemma3.1} (for finite $m>0$) or Lemma \ref{lemma3.2} (for $m=\infty$), let $H_i$ represent the mean curvature of a level set of $d_i$, $i\in \{ 1,2\}$, so that $\left ( H_X\right )_i=\Delta_X d_i$. Let $\delta=\Delta_X d_i \big\vert_{d_i=0}$. Then we see that $\left ( H_X\right )_i\le 0$ so long as $d_i$ is smooth, and holds in the barrier sense where it is not. Hence $\Delta_X (d_1+d_2)\le 0$. But $d_1+d_2$ will have an interior minimum and so must be constant on $D$. But since $\Delta_X d_i\le 0$ for $i\in \{ 1,2 \}$, then $\Delta_X d_i= 0$, and from \eqref{eq4.1} we deduce that $\tf \hess d_i=0$ and so $\hess d_i=0$. Then the level sets of the functions $d_i$ are totally geodesic and $U$ splits as claimed. Then from \eqref{eq4.1} we see that $\ric_X^m= 0$ and $X(d_i)=0$, so $X$ is tangent to the level sets $N$.
\end{proof}

Petersen and Wilhelm list the following cases.

\begin{corollary}(see \cite[Theorem 4]{PW})\label{corollary4.2}
Under the conditions of Lemma \ref{lemma4.1}, if $a>0$, then one of the following holds.
\begin{itemize}
\item [a)] Both $N_1$ and $N_2$ are $2$-sided in $M$ and $M\backslash {\bar U}$ consists of $2$ non-empty disjoint components, while ${\bar U}$ splits isometrically as a Riemannian product
\begin{equation}
\label{eq4.2}
{\bar U}\cong N_1\times [0,a].
\end{equation}
\item [b)] Both $N_1$ and $N_2$ are $2$-sided in $M$, $N_1$ and $N_2$ are isometric, and $M$ is isometric to a mapping torus
 \begin{equation}
\label{eq4.3}
M\cong\frac{N_1\times [0,a]}{(x,0)\sim (\phi(x),a)}
\end{equation}
where $\phi:N_1\to N_1$ is an isometry.
\item [c)] $N_1$ is $2$-sided but $N_2$ is $1$-sided and ${\bar U}$ splits isometrically as
\begin{equation}
\label{eq4.4}
{\bar U}\cong\frac{N_1\times [0,a]}{(x,0)\sim (\phi(x),a)}
\end{equation}
where now $\phi:N_1\to N_2$ is $2$-to-$1$ and a local isometry.
\item [d)] Both $N_1$ and $N_2$ are $1$-sided in $M$, and there are two $2$-to-$1$ local isometries $\phi_i:N\to N_i$ and $a_1<a_2$ such that $M$ splits isometrically as
\begin{equation}
\label{eq4.5}
M\cong\frac{N\times [a_1,a_2]}{(x,a_i)\sim (y,a_i)\text{ iff } \phi_1(x)=\phi_i(y)}\ .
\end{equation}
\end{itemize}
In every case, $X$ is tangent to the leaves (denoted by $N_1\times \{c\}$, $c\in [0,a]$, in (a--c) and by $N\times \{c\}$ in (d)).
\end{corollary}

\begin{proof}
Lemma \ref{lemma4.1} can be used to replace the first paragraph of the proof of \cite[Theorem 4]{PW}. The rest of the proof of \cite[Theorem 4]{PW} is essentially unchanged, and yields the theorem.
\end{proof}

Consider now a $2$-sided nonseparating totally geodesic surface $\Sigma$ in a closed manifold $M$ (\cite[Theorem 2.5]{CF}, \cite[Theorem (c)]{Galloway}). Choe and Fraser define a smooth function $M\backslash\Sigma \to {\mathbb R}$ that equals $0$ on a collar neighbourhood of one side of $\Sigma$ and $1$ on a collar of the other side. They then identify the range modulo the integers, obtaining a map which extends to all of $M$ smoothly, yielding a surjection $f:M\to S^1$. The induced map of fundamental groups $f_*:\pi_1(M)\to\pi_1(S^1)=Z$ counts the number of times a loop passes through $\Sigma$ from, say, left to right (minus the number of passages from right to left). The subgroup $G=\ker f_*$ consists of classes of loops which can be deformed so as never to pass through $\Sigma$. Then they consider the cover ${\hat M}:={\tilde M}/G$ where ${\tilde M}$ is the universal cover of $M$. The fundamental group of ${\hat M}$ consists of classes of those loops in $M$ that can be deformed to lie entirely in $\Sigma$.

\begin{lemma}\label{lemma4.3}
Let $M$ be a closed connected manifold with a vector field $X$ and an $m\in (0,\infty]$ such that $\ric_X^m \geq 0$ pointwise on $M$. Let $\Sigma$ be a closed, connected, embedded $2$-sided $X$-minimal hypersurface such that $M\backslash\Sigma$ is connected. Then $M$ is isometric to a mapping torus
 \begin{equation}
\label{eq4.6}
M\cong\frac{\Sigma \times [0,a]}{(x,0)\sim (\phi(x),a)}.
\end{equation}
\end{lemma}

\begin{proof}
As in \cite{CF}, construct the cover ${\hat M}$ of $M$. Since $\Sigma$ is $2$-sided and nonseparating, there must be a nontrivial loop, so the covering is nontrivial. The projection $\Pi :{\hat M}\to M$ induces local isometries, so $\ric_X^m \geq 0$ pointwise on ${\hat M}$. Consider two adjacent copies of $\Sigma$ in ${\mathbb M}$ and the closed region ${\bar U}$ between them. Then the argument proving Corollary \ref{corollary4.2}.(a) applies, so ${\bar U}\cong \Sigma\times [0,a]$ for some $a>0$. Furthermore, since ${\bar U}\times \{ 0 \}$ and ${\bar U}\times \{ a \}$ are both isometric to $\Sigma$, they are isometric to each other. Let $\phi$ denote this isometry. Identifying points via $\phi$ yields \eqref{eq4.6}.
\end{proof}

\begin{lemma}\label{lemma4.4}
Let $M$ be a closed connected manifold with a vector field $X$ and an $m\in (0,\infty]$ such that $\ric_X^m \geq 0$ pointwise on $M$. Let $\Sigma$ be a closed, connected, embedded $2$-sided $X$-minimal hypersurface such that $M\backslash\Sigma$ is disconnected. Let $D_1$, $D_2$ denote the disjoint connected components of $M\backslash \Sigma$. Let $i_j:\Sigma\to {\bar D}_j$, $I_j:D_j\to M$, and ${\mathcal I}:\Sigma\to M$ be the inclusions, where $j\in \{1,2\}$. Then the induced maps of fundamental groups $\left (i_j \right )_* :\pi_1(\Sigma)\to \pi_1({\bar D}_j)$, $\left (I_j\right )_*:\pi_1(D_j)\to M$, and $\left ({\mathcal I} \right )_* :\pi_1(\Sigma)\to M$ are all surjective.
\end{lemma}

\begin{proof}
The proof is the same as the proof of \cite[Theorem 2.5.(b)]{CF}. The only modification needed is that the superharmonicity condition $\Delta (d'+d'')\le 0$ quoted in that paper (for distance functions $d'$, $d''$ defined in that reference) becomes $X$-superharmonicity $\Delta_X (d'+d'')\le 0$, as used in the proof of our Lemma \ref{lemma4.1}.
\end{proof}

\begin{proof}[Proof of Theorem \ref{theorem1.4}] Theorem \ref{theorem1.4} is now an immediate corollary of Lemmata \ref{lemma4.3} and \ref{lemma4.4}.
\end{proof}

\section{Proofs of the splitting theorems}
\setcounter{equation}{0}

\subsection{Boundaries with more than one connected component: Theorems \ref{theorem1.5} and \ref{theorem1.6}}

\begin{proof}[Proof of Theorem \ref{theorem1.5}]
Let $D=d_M(N_1,N_2)$. Define distance functions $\rho_{N_i}:=d_M(p,N_i)$, $i\in \{1,2\}$. Let $F=\rho_{N_1}+\rho_{N_2}$. Let ${\tilde \Omega}$ be the subset of $\interior M$, the \emph{interior} of $M$, containing neither cut points nor shadow points for either $N_1$ or $N_2$ (a shadow point $x$ is one for which the minimizing curve from one boundary component to $x$ touches the other boundary component \emph{en route}). $F$ is smooth on ${\tilde \Omega}$, and $\Delta_X F = \Delta_X \rho_{N_1}+\Delta_X \rho_{N_2}\le 0$ by Lemma \ref{lemma3.1}, since $\Delta_X \rho_{N_i}$ are the $X$-mean curvatures $H_X$ of level sets of the functions $\rho_{N_i}$, $i\in \{1,2\}$. But $F$ attains a minimum $F=D$ at some $x\in M$ and then necessarily $F=D$ along the entire minimizing geodesic through $x$ from $N_1$ to $N_2$ (note that $x$ cannot be a shadow point). So $F$ attains its minimum at some $y\in {\tilde \Omega}\subset \interior M$. Then by the maximum principle $F$ is constant, and then $\Delta_X F = 0$ so $\Delta_X \rho_{N_1}=-\Delta_X \rho_{N_2}=(-1)^i\sqrt{(n-1)(n+m-1)}\delta$ for $i\in \{ 1,2 \}$.

We can now write $H_X:=\Delta_X\rho_{N_1}= -\sqrt{(n-1)(n+m-1)}\delta$ and invoke the equality statement in Lemma \ref{lemma3.1}. Then the tracefree part of $\hess \rho_{N_1}$ must vanish and $\ric_X^m(\nabla \rho_{N_1}, \nabla \rho_{N_1}) =-(n-1)\delta$. Since $\hess \rho_{N_1}$ is scalar, the level sets of $t$ are umbilic in $M$, and the twisted product splitting follows. The lemma also yields that $\sqrt {\frac{m}{n+m-1}}\Delta_X  \rho_{N_1} +\sqrt{\frac{n+m-1}{m}}X(\rho_{N_1}) =0$, so $g(X,\nabla \rho_{N_1})=m\sqrt{\frac{n-1}{n+m-1}}\delta$. Then $H:=\Delta \rho_{N_1}=\Delta_X \rho_{N_1} +g(X,\nabla \rho_{N_1}) =H_X+g(X,\nabla \rho_{N_1}) =-\sqrt{(n-1)(n+m-1)}\delta +m\sqrt{\frac{n-1}{n+m-1}}\delta =-\frac{(n-1)^{3/2}\delta}{\sqrt{n+m-1}} = -(n-1)c\delta$. Integrating, we obtain that the first fundamental form on level sets $t=\rho_{N_1}$ is $h(t)=e^{-2c t}g_1$. Then the metric on $M$ splits as a warped product $ds^2=dt^2+e^{-2c\delta t}g_1$,  $t\in [0,\ell]$.

By assumption, for any $k>0$ and any unit vector $v$ we have that $-(n-1)\delta\le \ric_X^m(kv+\epsilon \nabla \rho_{N_1},kv+\epsilon \nabla \rho_{N_1})$ for $\epsilon=\pm 1$. Expanding this and using that $\ric_X^m(\nabla \rho_{N_1},\nabla \rho_{N_1})=-(n-1)\delta$ we obtain that $0\le \pm 2k \ric_X^m(v,\nabla \rho_{N_1})+k^2 \ric_X^m(v,v)\le \pm 2k \ric_X^m(v,\nabla \rho_{N_1})+k^2 \left \vert \ric_X^m(v,v)\right \vert$. Dividing by one factor of $k>0$, we may now write this as $-k\vert \ric_X^m(v,v)\vert \le 2\ric_X^m(v,\nabla \rho_{N_1})\le k\vert \ric_X^m(v,v)\vert$. Then $k\searrow 0$ implies that $\ric_X^m(v,\nabla \rho_{N_1})=0$ for all $v$, so $\ric_X^m$ is block-diagonal and the condition $\ric_X^m$ descends to the restriction of $\ric_X^m$ to the orthogonal complement of $\frac{\partial}{\partial t}$. For $\delta=0$, we obtain $\ric_X^m(g_1)\ge 0$. For $\delta=1$, a brief calculation shows that the restriction of $\ric_X^m$  to $TN_1$ equals $\ric_{X^{\sharp}}^m(g_1)-(n-1)c^2 g_1$. Combining these cases and using $c^2=\frac{n-1}{n+m-1}$, we obtain that $\ric_{X^{\sharp}}^m(g_1)\ge -\frac{(n-1)^2}{(n+m-1)}\delta$.
\end{proof}

\begin{proof}[Proof of theorem \ref{theorem1.6}]
The first paragraph of the proof of Theorem \ref{theorem1.5} carries over to this situation, invoking Lemma \ref{lemma3.2} or Lemma \ref{lemma3.3} in place of Lemma \ref{lemma3.1} and concluding that $\Delta_X \rho_{N_1}=-\Delta_X \rho_{N_2}=0$. Then, as above, we have that $H_X=0$ at both boundaries and by the equality part of Lemma \ref{lemma3.2} or \ref{lemma3.3} we have that $\ric_X^m(\nabla_X \rho_{N_1},\nabla_X \rho_{N_1})=0$ and $\tf A=0$ where $A$ is the second fundamental form of level sets of $\rho_{N_1}$. Thus the level sets of $\rho_{N_1}$ are totally geodesic, so the metric splits as a product. Moreover, if $m\neq 1-n$, then by Lemma \ref{lemma3.3} we have $g(\nabla_X \rho_{N_1},X)=0$ and then $X$ must be tangent to the level sets of $\rho_{N_1}$.

Now for any $k>0$ and any unit vector $v$ we have that $0\le \ric_X^m(kv+\epsilon \nabla \rho_{N_1},kv+\epsilon \nabla \rho_{N_1})$ for $\epsilon=\pm 1$. Expanding this and using that $\ric_X^m(\nabla \rho_{N_1},\nabla \rho_{N_1})=0$ we obtain that $0\le \pm 2k \ric_X^m(v,\nabla \rho_{N_1})+k^2 \ric_X^m(v,v)\le \pm 2k \ric_X^m(v,\nabla \rho_{N_1})+k^2 \left \vert \ric_X^m(v,v)\right \vert$. As with the proof of Theorem \ref{theorem1.5}, we now divide by $k>0$, leaving one factor of $k$, to which we apply $k\searrow 0$ to conclude, as above, that the condition $\ric_X^m$ descends to the restriction of $\ric_X^m$ to the orthogonal complement of $\frac{\partial}{\partial t}$.
\end{proof}

\subsection{Complete manifolds with connected boundary: Theorems \ref{theorem1.7} and \ref{theorem1.8}}

\begin{proof}[Proof of Theorem \ref{theorem1.7}]
Let $\tau_p$ be the focal radius at $p\in N$, and let $\Omega:=\{ p\in N \big\vert \tau_p=\infty\}$. Because $M$ has an asymptotic end, $\Omega$ is non-empty. Because $\tau_p$ is continuous as a function of $p\in N$ and $N$ is compact, $\Omega$ is closed. But we will now show that $\Omega$ is open in $N$, and therefore $\Omega =N$. This in turn will imply that there are no focal points to $N$, so the normal exponential map generates a CMC foliation and so the metric splits as claimed.

Let $\gamma_{p_0}$ be a minimal geodesic issuing orthogonally from some point $p_0\in \Omega$. Then $\gamma_{p_0}$ is a ray from $N$. The distance function $d(p_0,\gamma_{p_0}(t))$ from $p_0$ will coincide with the distance from $N$. By Lemma \ref{lemma3.1} we have that the mean curvature of the distance function from $N$ obeys $H(t)\le -\sqrt{(n-1)(n+m-1)}\delta$. Therefore
\begin{equation}
\label{eq5.1}
\Delta_X d_1(q):=\Delta_X d(p_0,q)\le -\sqrt{(n-1)(n+m-1)}\delta
\end{equation}
at $q=\gamma_{p_0}(t)$, where the notation $d_1(q)=\inf_{\in N} \dist(x,q)=:\dist(N,q)$ follows the notation of \cite{CK}.

On the other hand, define the Busemann function at $q$ by $b^{\gamma_{p_0}}(q):=d_2(q):=\lim_{t\to\infty} \left \{ d(q,\gamma_{p_0}(t))-t\right \}$. The notation $d_2$ was introduced in \cite{CK}. By its definition, for any $x,y\in M$ we have $\left \vert b^{\gamma_{p_0}}(x) - b^{\gamma_{p_0}}(y)\right \vert \le d(x,y)$. (For this and other facts about Busemann functions, see for example \cite[p 286]{Petersen}.) The Busemann function thus defined may not be smooth near $q$, so we consider instead support functions defined as follows. Let $\eta_q:[0,\infty)\to M$ be an \emph{asymptote} for $\gamma_{p_0}$ with initial endpoint $q$; i.e., $\eta_q$ is an accumulation curve constructed by taking a sequence of curves joining $q$ to $\gamma_{p_0}(t)$ as $t\to\infty$. By construction, $b^{\gamma_{p_0}}(\eta_q(t)) =b^{\gamma_{p_0}}(\eta_q(0)) -t$ for all $t\in [0,\infty)$. Then we define
\begin{equation}
\label{eq5.2}
\begin{split}
b_t^{\gamma_{p_0}}(p):= &\, b^{\gamma_{p_0}}(p) -t +d(p,\eta_q(t))\\
\ge &\, b^{\gamma_{p_0}}(p) -t +\left \vert b^{\gamma_{p_0}}(q)- b^{\gamma_{p_0}}(\eta_q(t))\right \vert\\
\ge &\, b^{\gamma_{p_0}}(p) -t + b^{\gamma_{p_0}}(q)- b^{\gamma_{p_0}}(\eta_q(t))\\
= &\, b^{\gamma_{p_0}}(p) \equiv d_2(p),
\end{split}
\end{equation}
so that for $p$ near $q$, $b_t^{\gamma_{p_0}}$ is a smooth support function for $d_2$ from above.
Now we compute, using Lemma \ref{lemma3.6}, that
\begin{equation}
\label{eq5.3}
\Delta_X b_t^{\gamma_{p}} = \Delta_X d(p,\eta_q(t)) \le \sqrt{(n-1)(n+m-1)}y(t),
\end{equation}
where $y(t)$ is defined in equation \eqref{eq3.17}. We now define $F(p):= d_1(p)+d_2(p)$. Add \eqref{eq5.1} and \eqref{eq5.2} (at a common point $p$), and use that the definition of $y(t)$ implies that $\lim_{t\to\infty} y(t)= (n-1)\delta$. We obtain
\begin{equation}
\label{eq5.4}
\begin{split}
&\, \Delta_X \left (  d_1(p)+d_2(p) \right ) \le f(t)\to 0 \quad
\implies \quad \Delta_X \left (  d_1(p)+d_2(p)\right ) =\Delta_X F(p) \le 0,
\end{split}
\end{equation}
in the support sense on the interior of $M$ (precisely, on the open set excluding the cut locus of $N$, the ``shadow points'' of $N$, and $N$ itself; see \cite{CK}).

As in \cite[p 574]{CK}, since $N$ is compact, there is an $a_0\in N$ which minimizes $F\vert_N$. Let ${\mathcal F}$ denote the set of minima of $F$ in $M\backslash N$. If $x\in M$ and $y\in N$ such that $d(x,y)=d(x,N)=d_1(x)$, then $F(x)=d_1(x) +d_2(x)=d(x,y)+d_2(x)\ge d_2(y)\ge a_0$, with equality iff $y\in N\cap B_{a_0}(x)$ and $x$ lies on an asymptote that begins at $y$. Hence $x\in {\mathcal F}$ iff $x$ lies on an asymptote starting from $N$. Then by Calabi's maximum principle applied to \eqref{eq5.4} about $x\in {\mathcal F}$, there is an open set ${\mathcal O}\ni x$ such that ${\mathcal O}\subset {\mathcal F}$, so ${\mathcal F}$ is open in $M$. But then each point in ${\mathcal F}$ lies on an asymptote beginning on $N$, so the set of these initial endpoints is open as well.

Since ${\mathcal F}=N$, each point of $M=[0,\infty)\times N$ lies on a unique minimizing geodesic that leaves $N$ orthogonally. But since $N$ is compact, by Lemma \ref{lemma3.1}.(b) the focal radius of $N$ will be bounded above unless $H_X(0)=-\delta\sqrt{(n-1)(n+m-1)}$ at every point of $N$. But this must also be true at any $t_1>0$ as well, as we see by performing the translation $t\mapsto \tau=t-t_1$ and then invoking Lemma \ref{lemma3.1}.(b) at $\tau=0$. Hence we have $H_X(t)= -\delta \sqrt{(n-1)(n+m-1)}$ for all $t\in [0,\infty)$. Furthermore, Lemma \ref{lemma3.1} implies as well that we have $\ric_X^m(\gamma',\gamma')=-(n-1)\delta$, $\tf A \equiv 0$, and $g(\gamma',X)=-m\delta \sqrt{\frac{n-1}{n+m-1}}$. Writing $H_X(t)= - \sqrt{(n-1)(n+m-1)}\delta =H-g(\gamma',X)$, we find that $H= -\delta \sqrt{\frac{n-1}{n+m-1}}$, which together with $\tf A=0$ implies the product ($\delta=0$) or warped product ($\delta=1$) splitting
\begin{equation}
\label{eq5.5}
g=dt^2+e^{-2t\delta/\sqrt{(n-1)(n+m-1)}}g_N
\end{equation}
on $[0,\infty)\times N$, where $g_N$ is the induced metric on $t$-level sets diffeomorphic to $N$ and is independent of $t$.
\end{proof}

When $\delta=0$, if the metric splits as a product then the curvature condition $\ric_X^m(g)\ge 0$ implies $\ric_{X^\sharp}^m(g_N)\ge 0$. If $\delta=1$, a quick calculation using the warped product splitting yields $\ric(g_N)=\ric(g)\big\vert_{TN} +\frac{e^{-2t/\sqrt{(n-1)(n+m-1)}}g_N}{(n+m-1)}$. We can now also compute that $\pounds_{X^{\sharp}}g_N=\pounds_X g \big\vert_{TN}-\frac{2m e^{-2t/\sqrt{(n-1)(n+m-1)}}g_N}{(n+m-1)}$, where $X^{\sharp}$ is the restriction of $X$ to $TN$. Therefore when $\delta=1$ we have
\begin{equation}
\label{eq5.6}
\begin{split}
\ric_{X^{\sharp}}^m(g_N)\equiv &\, \ric(h)+\frac12 \pounds_{X^{\sharp}} h -\frac{1}{m}X^{\sharp}\otimes X^{\sharp} =\ric_X^m (g)\big\vert_{TN}- \frac{(m-1)}{(n+m-1)} e^{-2t/\sqrt{(n-1)(n+m-1)}}g_N\\
\ge &\, - \left ( n-1 +\frac{(m-1)}{(n+m-1)}\right ) e^{-2t/\sqrt{(n-1)(n+m-1)}}g_N\\
= &\, - \frac{n(n+m-2)}{(n+m-1)}e^{-2t/\sqrt{(n-1)(n+m-1)}}g_N \to 0 \text{ for }t\to\infty.
\end{split}
\end{equation}


\begin{proof}[Proof of Theorem \ref{theorem1.8}]
We follow the proof of Theorem \ref{theorem1.7}, replacing Lemma \ref{lemma3.1} in the second paragraph of that proof by Lemma \ref{lemma3.2} (if $m=\infty$) or \ref{lemma3.3} (if $m\in (-\infty,1-n]$), to ensure that $\Delta_X d(p_0,q)\le 0$. As well, Lemma \ref{lemma3.7} replaces Lemma \ref{lemma3.6} in the third paragraph of that proof, and then in place of equation \eqref{eq5.4} we have
\begin{equation}
\label{eq5.7}
\Delta_X \left (  d_1(p)+d_2(p)\right ) \le \frac{(n-1)e^{4K/(n-1)}}{t}\to 0
\end{equation}
as $t\to\infty$. We may then conclude, as in the proof of Theorem \ref{theorem1.7}, that ${\mathcal F}=N$; i.e., that the focal distance for geodesics leaving $N$ orthogonally is infinite at every point. Appealing again to Lemma \ref{lemma3.2} (if $m=\infty$) or Lemma \ref{lemma3.3} (if $m\in (-\infty,1-n]$), we then conclude that $H_X(0)=0$ and, as before, that in fact $H_X(t)=0$ for all $t$. But then, continuing to invoke these lemmata, we must have that $\ric_X^m(\gamma',\gamma')=0$, $\tf A=0$, and either $g(\gamma',X)=0$ along every geodesic $\gamma$ leaving $N$ orthogonally or $m=1-n$. If $m\neq 1-n$ then we obtain $H(t)=0$ and so the metric splits as a product.

If $m=1-n$, then we have that the second fundamental form is scalar (i.e., pure trace, or umbilic) with $H(t)=g(\gamma',X)$. This yields a twisted product splitting with the metric given in Gaussian normal coordinates relative to $N$ by \eqref{eq1.5}. Now it's elementary that if $B$ is a symmetric $(0,2)$-tensor obeying $B(w,w)\ge 0$ for all $w$ and $B(v,v)=0$ for a fixed nonzero vector $v$, then $B(v,w)=0$ for all $w$, and so we have that $\ric_X^{1-n}(\gamma',\cdot)= 0$. On the other hand, computing from \eqref{eq1.5} and writing the coordinates as $(t,x^i)$ where the $x^i$ are coordinates on $N$, we find that $\ric(\partial_t,\partial_i)=\frac{(2-n)}{(n-1)}\partial_i X_t$ and $\pounds_X g(\partial_t,\partial_i)=\partial_tX_i-\partial_iX_t-\frac{1}{(n-1)}X_tX_i$. This yields
\begin{equation}
\label{eq5.9}
\ric_X^m(\partial_t,\partial_i) =\frac12 \left [ \partial_tX_i-\frac{(n-3)}{(n-1)}\partial_iX_t\right ] -\frac{(n+m-1)}{m(n-1)}X_tX_i.
\end{equation}
Since $\ric_X^m(\partial_t,\partial_i)$  must vanish whenever $\ric_X^m\ge 0 $ and $\ric_X^m(\partial_t,\partial_t)=0$, so must the right-hand side of \eqref{eq5.9}, yielding $n-1$ differential equations on the $n$ components of $X$. Now set $m=1-n$. We obtain
\begin{equation}
\label{eq5.10}
0= \partial_tX_i-\frac{(n-3)}{(n-1)}\partial_iX_t.
\end{equation}
As an aside, the special case of $n=3$ is worthy of note. But for arbitrary $n>1$, now let $X$ be closed, so that in particular $\partial_tX_i=\partial_iX_t$. Then \eqref{eq5.10} reduces to $0=\partial_iX_t$ (compare \cite[Proposition 2.2]{Wylie} for the $X=df$ case), so $X_t\equiv g(\gamma',X)$ is independent of the coordinates $x^i$ on $N$, and so \eqref{eq1.5} is a warped product.
\end{proof}

\end{document}